\documentclass[12pt,reqno]{article}

\usepackage[usenames]{color}
\usepackage{amssymb}
\usepackage{graphicx}
\usepackage{amscd}

\usepackage[colorlinks=false]{hyperref}

\definecolor{webgreen}{rgb}{0,.5,0}
\definecolor{webbrown}{rgb}{.6,0,0}

\usepackage{color}
\usepackage{fullpage}
\usepackage{float}

\usepackage{graphics,amsmath,amssymb}
\usepackage{amsthm}
\usepackage{amsfonts}
\usepackage{latexsym}
\usepackage{epsf}

\newcommand{\seqnum}[1]{\href{http://oeis.org/#1}{\underline{#1}}}

\graphicspath{{converted_graphics/}}

\begin{document}

\begin{center}
\vskip 0.5cm
\end{center}

\theoremstyle{plain}
\newtheorem{theorem}{Theorem}
\newtheorem{corollary}[theorem]{Corollary}
\newtheorem{lemma}[theorem]{Lemma}
\newtheorem{proposition}[theorem]{Proposition}

\theoremstyle{definition}
\newtheorem{definition}[theorem]{Definition}
\newtheorem{example}[theorem]{Example}
\newtheorem{conjecture}[theorem]{Conjecture}

\theoremstyle{remark}
\newtheorem{remark}[theorem]{Remark}

\newcommand{\mean}{\mathop{\mathrm{mean}}}
\newcommand{\li}{\mathop{\mathrm{li}}}
\newcommand{\Exp}{\mathop{\mathrm{Exp}}}
\newcommand{\Gumbel}{\mathop{\mathrm{Gumbel}}}

\begin{center}
\vskip 2.2cm{\LARGE\bf 
On the Distribution of Maximal Gaps\\
\vskip .1in
Between Primes in Residue Classes
}
\vskip 0.7cm
\large
Alexei Kourbatov\\
www.JavaScripter.net/math\\
{\tt akourbatov@gmail.com}
\end{center}

\vskip 1cm

\begin{abstract} \noindent
Let $q>r\ge1$ be coprime positive integers.
We empirically study the maximal gaps $G_{q,r}(x)$ between primes $p=qn+r\le x$, $n\in{\mathbb N}$. 
Extensive computations suggest that almost always $G_{q,r}(x)<\varphi(q)\log^2x$.
More precisely, the vast majority of maximal gaps are near a trend curve $T$
predicted using a generalization of Wolf's conjecture:
$$
G_{q,r}(x) ~\sim~ T(q,x) = {\varphi(q)x\over\li x} 
\Big( 2\log{\li x\over\varphi(q)} - \log x + b\Big),
$$ 
where $b=b(q,x)=O_q(1)$.
The distribution of properly rescaled maximal gaps $G_{q,r}(x)$ is close to 
the Gumbel extreme value distribution. However, the question whether 
there exists a limiting distribution of $G_{q,r}(x)$ is open. 
We discuss possible generalizations of Cram\'er's, Shanks, and Firoozbakht's
conjectures to primes in residue classes.
\end{abstract}

\vskip 0.7cm


\pagebreak
\section{Introduction}

Let $q$ and $r$ be positive integers such that $1\le r < q$ and $\gcd(q,r)=1$.
A {\em residue class} is a set of integers that leave the same remainder $r$ 
when divided by a given modulus $q$.
Dirichlet proved in 1837 that for natural $n$ the integer $qn+r$ 
is prime infinitely often; this is Dirichlet's theorem on~arithmetic progressions.
The prime number theorem tells us that the total number of primes $\le x$ 
is asymptotic to the logarithmic integral of $x$, 
i.e., $\pi(x)\sim\li x$. Moreover,
the generalized Riemann hypothesis implies that primes are distributed approximately equally 
among the $\varphi(q)$ residue classes modulo $q$ corresponding to specific values of $r$. 
Thus each residue class contains a positive proportion, about $1\over\varphi(q)$, of all primes below $x$.
Accordingly, the average prime gap below $x$ is about $x/\li x \sim \log x$,
while the average gap between primes $p=qn+r\le x$ is about $\varphi(q)x/\li x \sim \varphi(q)\log x$.

In this paper we empirically study the growth and distribution of the values of the function 
$G_{q,r}(x)$, the record (maximal) gap between primes of the form $qn+r$ below $x$, for $x<10^{12}$.
The special case $G_{2,1}(x)$, i.e.\ the maximal gaps between primes below $x$, 
has been studied by many authors; 
see e.\,g.~\cite{bhp2001,cadwell,cramer,fgkmt,nicely,toes2014,shanks,wolf2011,wolf2014}.
Computational experiments investigating the actual distribution of (properly rescaled) values $G_{q,r}(x)$
are of interest, in part, because in Cram\'er's {\em probabilistic model of primes} \cite{cramer}
there exists a limiting distribution of maximal prime gaps \cite{kourbatov2014}, 
namely, the Gumbel extreme value distribution.

The fact that values of certain arithmetic functions have limiting distributions
is among the most beautiful results in number theory.
A well-known example is the limiting distribution of $\omega(n)$, the number of distinct prime factors of $n$.
The {\em Erd\H{o}s-Kac theorem} states that, roughly speaking, the values of $\omega(n)$ for $n\le x$  
follow the {\em normal distribution} with the mean $\log\log x$ and standard deviation $\sqrt{\log\log x}$, as $x\to\infty$.
Note that it is virtually impossible to observe the Erd\H{o}s-Kac normal curve for $\omega(n)$ 
in a computational experiment because such an experiment would involve factoring a lot of gigantic values of $n$.
By contrast, computations described in Section \ref{numresults} and {\it Appendix} are quite manageable 
and do allow one to draw maximal gap histograms, which turn out to closely fit the Gumbel distribution.
These results, together with \cite{kourbatov2014}, support the hypothesis 
that there is a limiting distribution of (properly rescaled) values of $G_{q,r}(x)$.
However, a formal proof or disproof of existence of a Gumbel limit law for 
maximal gaps between primes in residue classes seems beyond reach.

\section{Notation and abbreviations}

\begin{tabular}{ll}  
$p_k$             & the $k$-th prime;\, $\{p_k\} = \{2,3,5,7,11,\ldots\}$                                  \\
$\pi(x)$          & the prime counting function: the total number of primes $p_k\le x$                     \\
$\pi_{q,r}(x)$    & the prime counting function in residue class $r$ modulo $q$:                           \\
                  & the total number of primes $p=qn+r\le x$, $n\in{\mathbb N}^0$                          \\
$\gcd(q,r)$       & the greatest common divisor of $q$ and $r$                                             \\
$\varphi(n)$      & Euler's $\varphi$ function: the number of positive $m\le n$ with $\gcd(m,n)=1$         \\
$G(x)$            & the record (maximal) gap between primes $\le x$                                        \\
$G_{q,r}(x)$      & the record (maximal) gap between primes $p=qn+r \le x$                                 \\
$a=a(q,x)$        & the expected average gap between primes $p=qn+r\le x$;                                 \\
                  & defined as $a=\varphi(q)x/\li x$         \\ 
$b=b(q,x)$        & the correction term in equation (\ref{trend})                                          \\
$N(x)$            & the total number of maximal prime gaps $G$ with endpoints $p\le x$                     \\
$N_{q,r}(x)$      & the total number of maximal gaps $G_{q,r}$ with endpoints $p\le x$                     \\
i.i.d.            & independent and identically distributed                                                \\
cdf               & cumulative distribution function                                                       \\
pdf               & probability density function                                                           \\
EVT               & extreme value theory                                                                   \\
GRH               & generalized Riemann hypothesis                                                         \\
$\Exp(x;\alpha)$  & the exponential distribution cdf: \ $\Exp(x;\alpha)=1-e^{-x/\alpha}$                   \\
$\Gumbel(x;\alpha,\mu)$& the Gumbel distribution cdf: \ $\Gumbel(x;\alpha,\mu) =  
                    e^{-e^{-{x-\mu\over\vphantom{f}\alpha}}}$                                              \\
$\alpha$          & the {\em scale parameter} of exponential/Gumbel distributions, as applicable           \\
$\mu$             & the {\em location parameter} ({\em mode}) of the Gumbel distribution                   \\
$\gamma$          & the Euler-Mascheroni constant: \ $\gamma = 0.57721\ldots$                              \\
$\Pi_2$           & the twin prime constant: \ $\Pi_2 = 0.66016\ldots$                                     \\ 
$\log x$          & the natural logarithm of $x$                                                           \\
$\li x$           & the logarithmic integral of $x$: \
                    $\displaystyle\li x \,= \int_0^x{\negthinspace}{dt\over\log t}
                                        \,= \int_2^x{\negthinspace}{dt\over\log t} + 1.04516\ldots$        \\
$\zeta(s)$        & the Riemann $\zeta$-function: \ $\zeta(s)=\sum\limits_{k=1}^\infty{1\over k^s}$  \\
\end{tabular}

\pagebreak

\section{Numerical results}\label{numresults}

Using the PARI/GP program {\tt maxgap.gp} (see {\em Appendix}) we have computed the maximal gaps $G_{q,r}(x)$
for $x<10^{12}$, with many values of $q\in[4,10^5]$.
We used all admissible values of $r\in[1,q-1]$, $\gcd(q,r)=1$, 
to assemble a complete data set of maximal gaps for a given $q$.
This section summarizes our numerical results.

\medskip

\begin{center}Table~1. \ Example: record gaps between primes $p=1000n+1$ \\[0.5em]
\begin{tabular}{rrrr}
\hline {\large $\vphantom{1^{1^1}}$}
Start of gap  &  End of gap ($p$) &  Gap $G_{q,r}(p)$  & Rescaled gap $w$, eq.\,(\ref{rescaledw}) \\
[0.5ex]\hline
\vphantom{\fbox{$1^1$}}
         3001 &          4001 &     1000 & $-1.2053923292 $\phantom{111}\\
         4001 &          7001 &     3000 & $-0.7312116351 $\phantom{111}\\
         9001 &         13001 &     4000 & $-0.7329842896 $\phantom{111}\\
        28001 &         51001 &    23000 & $ 3.3272985119 $\phantom{111}\\
       294001 &        318001 &    24000 & $ 1.4713160904 $\phantom{111}\\
       607001 &        633001 &    26000 & $ 1.1176630346 $\phantom{111}\\
      4105001 &       4132001 &    27000 & $-0.7751904271 $\phantom{111}\\
      5316001 &       5352001 &    36000 & $ 0.5041782491 $\phantom{111}\\
     14383001 &       14424001&    41000 & $ 0.1458184959 $\phantom{111}\\
     26119001 &       26163001&    44000 & $-0.1041704135 $\phantom{111}\\
     46291001 &       46336001&    45000 & $-0.6486127233 $\phantom{111}\\
     70963001 &       71011001&    48000 & $-0.7238822913 $\phantom{111}\\
     95466001 &       95515001&    49000 & $-0.9421573040 $\phantom{111}\\
    114949001 &      115003001&    54000 & $-0.4549472918 $\phantom{111}\\
    229690001 &      229752001&    62000 & $-0.2205245619 $\phantom{111}\\
    242577001 &      242655001&    78000 & $ 1.9009667495 $\phantom{111}\\
    821872001 &      821958001&    86000 & $ 1.2346880297 $\phantom{111}\\
   3242455001 &     3242545001&    90000 & $-0.1834023117 $\phantom{111}\\
   7270461001 &     7270567001&   106000 & $ 0.5642767037 $\phantom{111}\\
  11281191001 &    11281302001&   111000 & $ 0.5051584212 $\phantom{111}\\
  32970586001 &    32970700001&   114000 & $-0.6900065927 $\phantom{111}\\
  50299917001 &    50300032001&   115000 & $-1.1757092375 $\phantom{111}\\
  63937221001 &    63937353001&   132000 & $ 0.2732357356 $\phantom{111}\\
  92751774001 &    92751909001&   135000 & $ 0.0415120508 $\phantom{111}\\
 286086588001 &   286086729001&   141000 & $-0.9866954146 $\phantom{111}\\
 334219620001 &   334219767001&   147000 & $-0.6218345554 $\phantom{111}\\
 554981875001 &   554982043001&   168000 & $ 0.6629883248 $\phantom{111}\\
1322542861001 &  1322543032001&   171000 & $-0.3590754946 $\phantom{111}\\
2599523890001 &  2599524073001&   183000 & $-0.2691207470 $\phantom{111}\\
4651789531001 &  4651789729001&   198000 & $ 0.1945079654 $\phantom{111}\\
7874438322001 &  7874438536001&   214000 & $ 0.7901492021 $\phantom{111}\\
8761032430001 &  8761032657001&   227000 & $ 1.7539727990 $\phantom{111}\\
\hline
\end{tabular}
\end{center}

\subsection{The growth trend of maximal gaps $G_{q,r}(x)$}

Let us begin with a simple example. For $q=1000$, $r=1$, running the program {\tt maxgap.gp} 
produces the results shown in Table 1. 
It is easy to check that all gaps $G_{q,r}(p)$ in the table satisfy the inequality
$G_{q,r}(p)<\varphi(q)\log^2 p$, suggesting several possible generalizations of 
{\it Cram\'er's conjecture} (see sect.\,\ref{gencram}).

\begin{figure}[h] 
  \centering
  \includegraphics[bb=4 2 420 280,width=5.6in,height=3.7in,keepaspectratio]{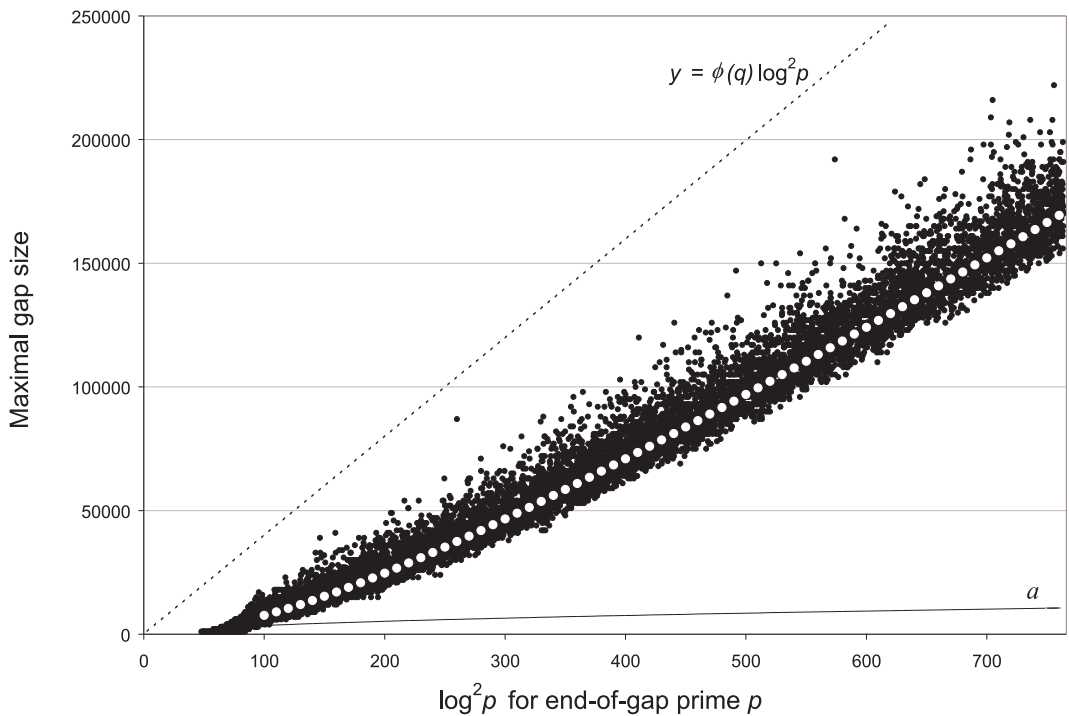}
  \caption{Record gaps $G_{q,r}(p)$ between primes $p=qn+r<10^{12}$ in residue classes mod $q$,
   with $q=1000$, $\varphi(q)=400$, $\gcd(q,r)=1$.  Plotted (bottom to top): 
   {\em average gaps} $a={\varphi(q)\cdot p\over\li p}$ between primes $\le p$ in residue classes mod $q$ (solid curve);
   {\em trend curve} $T$ of eq.\,(\ref{trend}) (white dotted curve);
   the conjectural (a.s.) upper bound for $G_{q,r}(p)$: $y=\varphi(q)\log^2p$ (dashed line).
  }
  \label{fig:1}
\end{figure}

Figure \ref{fig:1} shows all gaps $G_{q,r}(p)$ for $q=1000$, $\forall r\in[1,q-1]$, $\gcd(q,r)=1$, $p<10^{12}$. 
The horizontal axis $\log^2p$ reflects the actual end-of-gap prime $p=qn+r$ of each maximal gap.
The results for other values of $q$ closely resemble Fig.\,\ref{fig:1}.
We once again observe that $G_{q,r}(p)<\varphi(q)\log^2 p$. 
(For rare examples with $G_{q,r}(p)>\varphi(q)\log^2 p$, see Appendix \ref{appendix-extra-large}.)
\vskip 1mm

\noindent
The vast majority of record gaps $G_{q,r}(x)$ are near a smooth trend curve $T$:
\begin{equation}\label{trend}
G_{q,r}(x) ~\sim~ T(q,x) = a\cdot \left( 2\log{\li x\over\varphi(q)} - \log x + b\right),
\qquad\mbox{ cf.\,\cite{wolf2016},}
\end{equation}
where $a$ is the {\em expected average gap} between primes in a residue class mod $q$, defined as
\begin{equation}\label{defa}
a ~=~ a(q,x) ~=~ \varphi(q){x\over\li x}, 
\end{equation}
and $b=b(q,x)=O_q(1)$ is a correction term. Clearly, for $q=\mbox{const}\,$ and $x\to\infty$ we have
$$
T(q,x) ~\lesssim~ \varphi(q)\log^2 x \qquad\mbox{ and }\qquad
a(q,x) ~\lesssim~ \varphi(q)\log x. 
$$

The trend equation (\ref{trend}) is derived from a new generalization \cite{wolf2016} 
of Wolf's conjecture on maximal prime gaps \cite[Conjecture 3]{wolf1998}; 
see also \cite[eq.\,(31)]{wolf2011} and \cite[eq.\,(15)]{wolf2014}.
Wolf's heuristic implies that in eq.\,(\ref{trend}) the term $b$ tends to a constant for large $x$. 
Indeed, empirically we find
\begin{equation}\label{defb}
b ~=~ b(q,x) ~\approx\, \left(b_0+{b_1\over(\log\log x)^d}\right)\log\varphi(q),
\end{equation}
where the parameter values 
\begin{equation}\label{defb0}
b_0=1, \qquad b_1=4, \qquad d=2.7
\end{equation}
are close to optimal for $q\in[10^2,10^5]$ and $x\in[10^6,10^{12}]$. 
(Sometimes the parameter values may need a little adjustment, especially when $\varphi(q)$ is small;
usually it is enough to adjust $b_1$ and/or $d$ while keeping $b_0=1$.) 
We stress that eqs.\,(\ref{defb}),\,(\ref{defb0}) are purely empirical.

\subsection{The distribution of maximal gaps}\label{distsect}

We have thus found that the maximal gaps between primes in each residue class are mainly observed within
a strip of increasing width $O(a)$ around the trend curve $T(q,x)$ of eq.\,(\ref{trend}),
where $a=a(q,x)$ is the expected average gap between primes in the respective residue class.
Now let us take a closer look at the distribution of maximal gaps in
the neighborhood of this trend curve. We perform a rescaling transformation:
subtract the trend $T(q,x)$ from actual gap sizes, and then divide the result by the ``natural unit'' $a$.
All record gaps $G_{q,r}(p)$ are mapped to {\it rescaled values} $w$:
\begin{equation}\label{rescaledw}
G_{q,r}(x) ~~\mapsto~~ w = {G_{q,r}(x)-T(q,x)\over a(q,x)},
\qquad\mbox{ where } \ a(q,x) = {\varphi(q)x\over\li x}.
\end{equation} 

\begin{figure}[p] 
  \centering
  \includegraphics[bb=90 0 528 777,width=5.67in,height=7.8in,keepaspectratio]{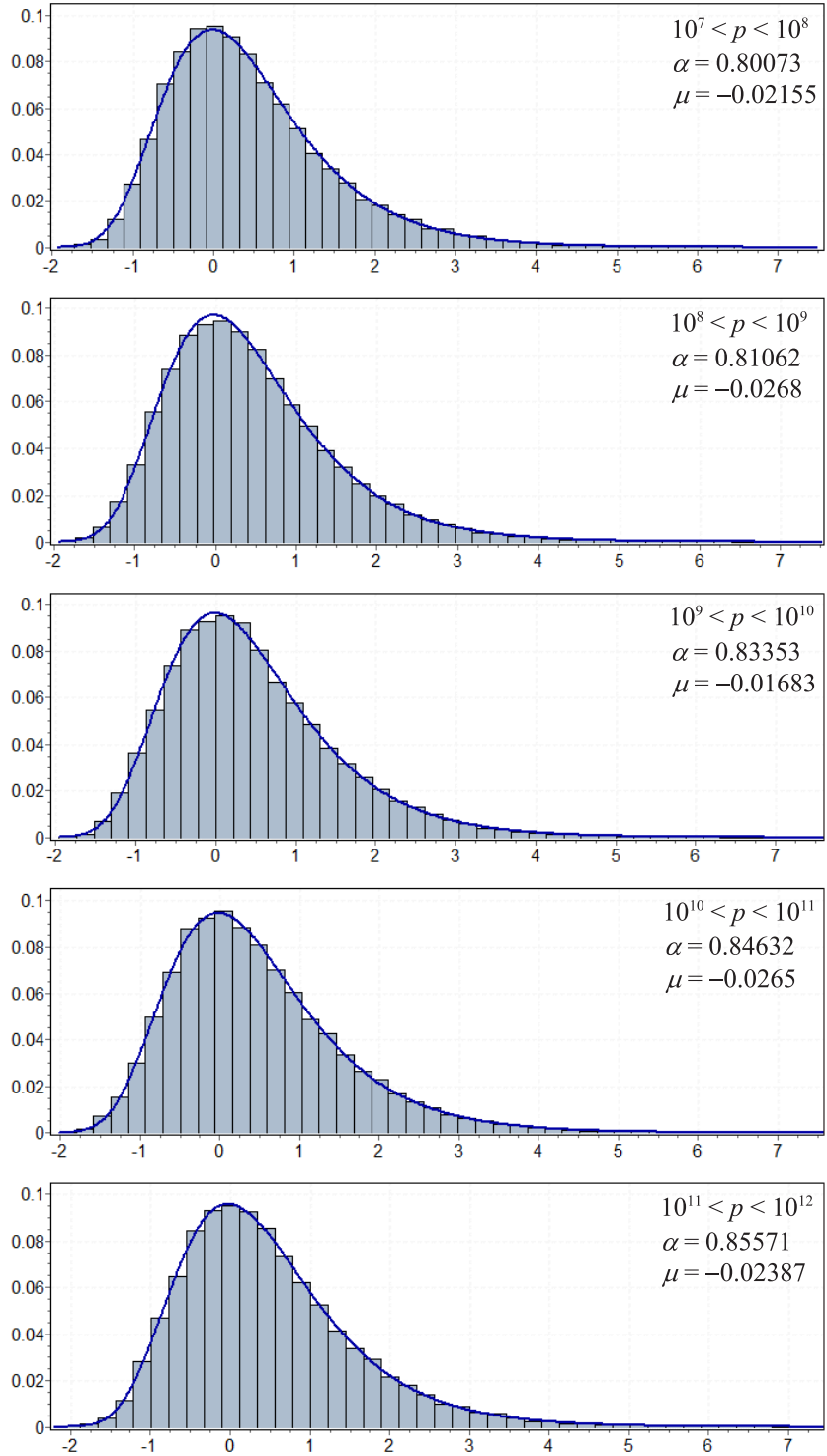}
  \caption{Histograms of rescaled maximal gaps $w$ and fitted Gumbel distributions (pdf)
   for 8-, 9-, 10-, 11-, and 12-digit primes $p=qn+r$, with $q=10007$, $r\in[1,10006]$. 
  }
  \label{fig:2}
\end{figure}

Figure \ref{fig:2} shows the histograms of rescaled maximal gaps for $q=10007$.
(Histograms for other $q$ look similar to Fig.\,\ref{fig:2}.)
We can see at once that the histograms and fitting distributions are skewed: 
the right tail is longer and heavier. This skewness is a well-known feature of extreme value distributions.
Among all two-parameter distributions supplied by the distribution fitting software \cite{easyfit},
the best fit is the {\em Gumbel distribution}. 
This opens up the question whether the Gumbel distribution is the limit law for properly rescaled
sequences of the $G_{q,r}(x)$ values as $x\to\infty$; cf.\,\cite{kourbatovOEIS50,kourbatov2014}. 
Does such a limiting distribution exist at all?

If we look at three-parameter distributions, then one of the best fits is the 
{\em Generalized Extreme Value} (GEV) distribution, which includes the Gumbel distribution as a special case.
The {\em shape parameter} in the best-fit GEV distributions is very close to zero; 
note that the Gumbel distribution is a GEV distribution whose shape parameter is exactly zero.
Thus, convergence to the Gumbel limit law appears feasible for
properly rescaled infinite sequences $G_{q,r}(p)$ as $p\to\infty$.
Clearly, a more precise knowledge of the trend function $T(q,x)$ in (\ref{trend})
is crucial for proper rescaling.

For best-fit Gumbel distributions computed for $10^k<p<10^{k+1}$, $7\le k\le 11$, with various values of $q<10^5$,
the Kolmogorov--Smirnov goodness-of-fit statistic slowly improves
from about 0.1 (when $\varphi(q)\approx 30$)
to about 0.03 (when $\varphi(q)\approx 300$),
to about 0.01 (when $\varphi(q)\approx 3000$),
to about 0.003 (when $\varphi(q)\approx 30000$).

There are notable differences between our situation shown in Fig.\,\ref{fig:2} 
and the distributions of maximal gaps between {\em prime $k$-tuples} shown in \cite[Fig.\,4]{kourbatov2013}~---
even though in either case the Gumbel distribution is a good fit, with the vast majority of record gaps
occurring within $\pm2a$ of the respective trend curve.
\begin{itemize}
\item In our case (Figs.\,\ref{fig:1} and \ref{fig:2}), there is plenty of data:
  many thousands of data points are available. 
  On the other hand, data on maximal gaps between prime $k$-tuples are scarce; 
  at present we have fewer than 100 data points for any given type of prime $k$-tuple 
  \cite{kourbatov2013}. 

\item In case of  gaps between prime $k$-tuples, $k\ge2$, 
  maximal gap sizes are usually not too far from $a\log(x/a)$,
  where $a$ denotes the expected average gap between $k$-tuples of a particular type.
  Thus a trend equation based on extreme value theory produces satisfactory results
  for maximal gaps between $k$-tuples; see \cite[sect.\,4-5]{kourbatov2013}.
  (In the {\tt randomgap.gp} control experiments, EVT also works well; see sect.\,\ref{controlexp}.)
  By contrast, for maximal gaps $G_{q,r}(p)$, eq.\,(\ref{trend}) appears to be a better fit to the observed 
  trend.\footnote{
   If, instead of using Wolf's conjecture-based equation (\ref{trend}), we were to model the trend of $G_{q,r}(p)$ 
   with an EVT-based equation akin to (\ref{rtrend}): 
   $G_{q,r}(x) \sim a\cdot(\log{\li x\over\varphi(q)}+\beta)$,
   then our trend equation would need an {\it unbounded} correction term $\beta=O(\log\log x)$;\,
   see also section \ref{genwolf}.
  } 
\end{itemize}

As noted by Brent \cite{brent2014}, primes seem to be less random than twin primes. 
We can add that, likewise, record gaps between primes in a residue class 
seem to be somewhat less random than those for prime $k$-tuples. Our analysis also shows that
primes $p=qn+r$ do not go quite as far from each other as in the {\tt randomgap.gp} model;
see section \ref{controlexp}.
Pintz \cite{pintz2007} discusses various other aspects of the ``random'' and not-so-random behavior of primes.

\subsection{Control experiments: records among random gaps}\label{controlexp}

For comparison, we performed a series of computational experiments studying
the records in a growing sequence of integers $p$ separated by {\em random gaps} 
$\lceil\xi\rceil$, where $\xi$ is an exponentially distributed random variable with mean $\varphi(q)\log p$.
In these experiments we used the PARI/GP program {\tt randomgap.gp} (see {\it Appendix}).

\medskip\noindent
{\bf Trend equation for record random gaps.}
We observed that sizes of maximal random gaps below $x$ mostly follow this trend:
\begin{equation}\label{rtrend}
\mbox{maximal random gap size } g ~\sim~
T_{\mbox{\tiny{rand}}}(q,x) 
      = {\varphi(q)x\over\li x} \cdot \log{\li x\over\varphi(q)}.
\end{equation}
  We can heuristically derive equation (\ref{rtrend}) as follows. 
  {\em Extreme value theory} predicts that, for $N$ consecutive events occurring at i.i.d.\ random intervals
  with cdf~$\Exp(\xi;\alpha)=1-e^{-\xi/\alpha}$ (i.\,e.\ at average intervals $\alpha$), 
  the most probable maximal interval between events is about $\alpha\log N$,
  while the width of extreme value distribution is $O(\alpha)$; see \cite{gumbel1958}. 
  As $N\to\infty$, EVT predicts that the distribution of maximal intervals approaches the 
  Gumbel distribution $\Gumbel(x;\alpha,\mu)$ with scale $\alpha$ and mode $\mu=\alpha\log N$.
  To derive eq.\,(\ref{rtrend}) we simply take the estimate $\alpha\log N$ with  
  $N\approx\li x/\varphi(q)$ and $\alpha\approx a(q,x)\approx x/N\approx \varphi(q)x/\li x$.
For $q=2$, expression (\ref{rtrend}) reduces to approximately $\log(x)\log(x/\log x)$; 
cf.~Cadwell \cite[p.\,912]{cadwell}.

\begin{figure}[h] 
  \centering
  \includegraphics[bb=2 0 642 171,width=6.2in,height=1.7in]{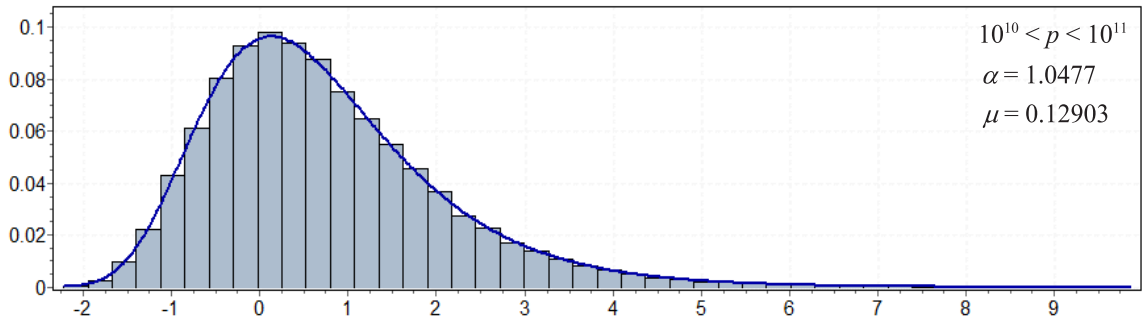}
  \caption{Histogram of rescaled values $h$ for records among exponentially distributed 
   {\em random} gaps obtained in a computational experiment with {\tt randomgap.gp}; $q=10007$. 
   The smooth curve is the Gumbel distribution (pdf) with scale $\alpha=1.0477$ and mode $\mu=0.12903$.}
  \label{fig:rand10007}
\end{figure}

\medskip\noindent
{\bf Rescaling transformation. }
The observed trend curve (\ref{rtrend}) for random gaps differs from $T(q,x)$ of eq.\,(\ref{trend}).
Accordingly, instead of eq.\,(\ref{rescaledw}), 
for maximal random gaps $g$ we have to use a different rescaling transformation:
\begin{equation}\label{rescaledh}
g ~~\mapsto~~ h 
             = {g - a\log{\li x\over\varphi(q)} \over a},
\end{equation}
where $a = a(q,x) = {\varphi(q)x\over\li x}$.
Figure \ref{fig:rand10007} shows a histogram of the rescaled values $h$ obtained for $q=10007$.
The histogram closely approximates the Gumbel distribution.

\medskip\noindent
{\bf Mode parameter $\mu$. }
It was observed that the mode parameter $\mu$ of the best-fit Gumbel distribution for 
the rescaled values $h$ was close to 0; in most computational experiments we obtained 
$|\mu|<0.3$. This means that the most probable values of record random gaps $g$ 
are $a\log{\li x\over\varphi(q)} + O(a)$ (before rescaling).

\medskip\noindent
{\bf Scale parameter $\alpha$. }
In the best-fit Gumbel distributions for the rescaled values $h$,
the scale parameter $\alpha$ was close to 1; in most experiments with random gaps we obtained 
$\alpha$ fluctuating in the range $0.9<\alpha<1.1$; see Figure~\ref{fig:rand10007}. 
On the other hand, our experiments of section \ref{distsect} 
produced best-fit Gumbel distributions for $w$ (the rescaled values of $G_{q,r}(x)$)
whose scale parameter $\alpha$ seemed to slowly approach 1 {\em from below} 
and was mostly in the range $0.7<\alpha<1$; see Figure~\ref{fig:2}.

\subsection{How many record gaps are there?}\label{howmanyrecords}

Let $N(x)$ be the number of maximal gaps observed between primes $\le x$.
Numerical data on the known maximal prime gaps allow us to make a conjecture:
below the $k$-th prime, there are about $2\log k$ maximal gaps, i.e., about twice as many as
the expected number of records in a sequence of $k$ i.i.d.\ random variables \cite[A005250, A005669]{oeis}.
Taking into account that $\pi(x)\sim\li x$, we can write
the above conjecture as
\begin{equation}\label{twologlix}
N(x) ~\sim~ 2\log\li x.
\end{equation}
Can we generalize this formula to maximal gaps between primes in residue class $r$ mod $q$?

Denote by $N_{q,r}(x)$ the total number of record gaps $G_{q,r}$ with endpoints below $x$.
Given a coprime pair $(q,r)$, we observe that usually
the interval $[x,ex]$ contains an endpoint of at least one record gap $G_{q,r}(p)$. 
(Of course, some intervals $[x,ex]$ contain no record gaps $G_{q,r}(p)$ for a given pair $(q,r)$.
For example, for $q=1000$ and $r=1$, the interval $[10^5,10^5e]$ does not contain endpoints 
of any record gaps $G_{1000,1}$; see Table 1.)

\vskip 6mm

\begin{figure}[h] 
  \centering
  \includegraphics[bb=6 0 430 255,width=5.7in,height=3.4in]{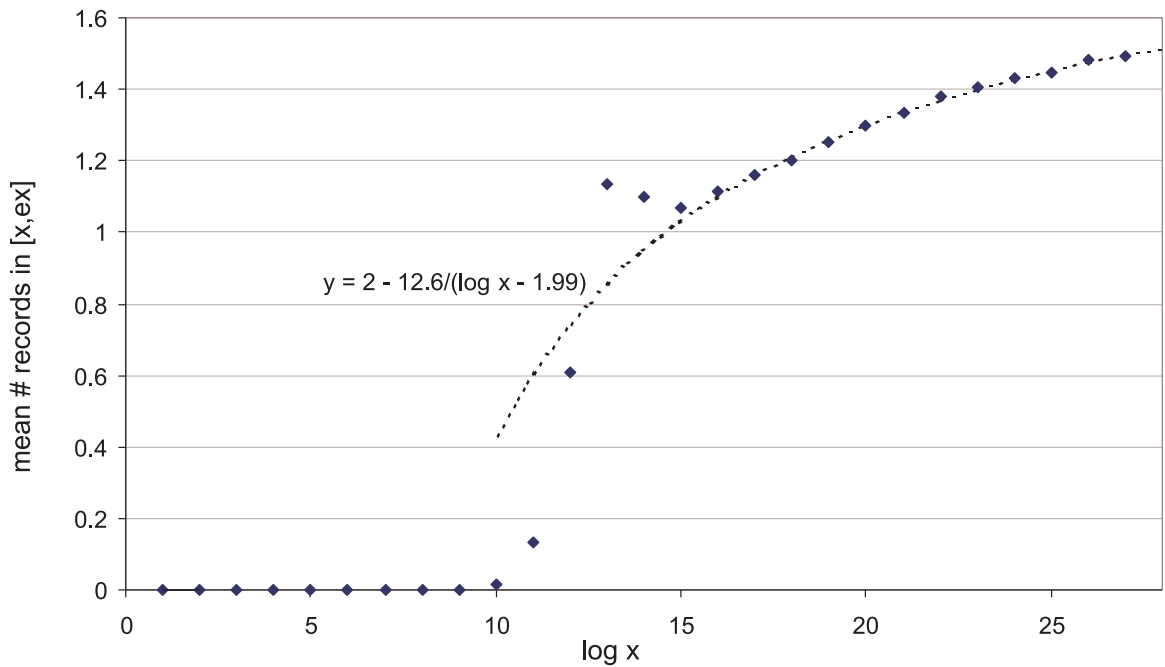}
  \caption{The average numbers of record gaps $G_{q,r}(p)$ with end-of-gap primes $p\in[e^k,e^{k+1}]$
           computed for $q=20011$, $k\le27$. Averaging was performed for all $r\in[1,20010]$.}
  \label{fig:4}
\end{figure}

A closer examination reveals that the average number of records in the interval $[x,ex]$
slowly grows as $x\to\infty$. To quantify this phenomenon, we used a modified version
of the program {\tt maxgap.gp} to compute average numbers of record gaps with endpoints $p\in[e^k,e^{k+1}]$ for 
many different $q$. With each $q$, the averaging was performed for all $r\in[1,q-1]$, $\gcd(q,r)=1$.
Figure \ref{fig:4} gives the results of this computation for $q=20011$.

Analysis of computation results shows that the average number of records in the interval $[x,ex]$ 
is satisfactorily approximated by this rational function of $\log x$ ({\it hyperbola} in Fig.\,\ref{fig:4}):
\begin{equation}
\mean_r \big( N_{q,r}(ex) - N_{q,r}(x) \big)~\approx~ 2 - {\kappa(q)\over\log x - \delta(q)}.
\end{equation}
The leading term 2 in this approximation is heuristically justified \cite[sect.\,2.3]{kourbatov2017}.
But if we want to construct an empirical (almost sure) {\em upper bound} for 
the number of records in the interval $[x,ex]$, we can simply use a function of the form 
$(\log \log x)^\lambda$ with some $\lambda>0$.
We will take advantage of this observation in the next section; 
see conjectures \ref{conjnum1} and \ref{conjnum2} on the number of record gaps.
If, in addition, we can assume the existence of the limiting distribution of $G_{q,r}(x)$, 
this will allow us to predict how often $G_{q,r}(x)>\varphi(q)\log^2x$; 
see section \ref{consequences}.

\section{New conjectures}

\subsection{Conjectures on the number of record (maximal) gaps}\label{conjnum}

It is well known that in a sequence of $n$ i.i.d.\ random variables with a continuous pdf
the expected number of record values grows asymptotically as fast as $\log n$; see e.g. \cite{krug2007}.
Indeed, 
\begin{itemize}
 \item the 1st term in an i.i.d.\ random sequence is guaranteed to be a record,
 \item the 2nd term is a record with probability ${1\over2},\ldots$
 \item the $k$th term is a record with probability ${1\over k}$ (by symmetry).
\end{itemize}
Thus the expected total number of records in an i.i.d.\ random sequence of $n$ elements is
$$
1 + {1\over2} + {1\over3} + \ldots + {1\over n}
~=~ \sum\limits_{k=1}^n {1\over k} ~\sim~ \log n + \gamma.
$$

As before, let $N_{q,r}(x)$ be the total number of record (maximal) gaps observed below $x$ in the sequence $G_{q,r}$.
Although the gaps between primes of the form $p=qn+r$ are {\bf\em not random variables}, computations show that 
the total numbers of records $N_{q,r}(x)$ behave somewhat like those for i.i.d.\ random variables.
Based on computational results we make the following conjectures.

\subsubsection{Conjectures on the total number of records $N_{q,r}(x)$}\label{conjnum1} 

{\bf Weak conjecture.}
For any $q>r\ge1$ with $\gcd(q,r)=1$, there exists $\lambda>0$ such that
$$
N_{q,r}(x) < \log x (\log\log x)^\lambda \qquad\mbox{almost always. }
$$

\smallskip\noindent
{\bf Strong conjecture.}
For any $q>r\ge1$ with $\gcd(q,r)=1$, there exists $C>2$ such that
$$
N_{q,r}(x) < C\log x.
$$
{\it Remark.} Here is a straightforward generalization of ({\ref{twologlix}) for ``guesstimating'' $N_{q,r}(x)$:
\begin{equation}\label{twologlixphiq}
N_{q,r}(x) ~\approx~ \max\left(0,\,2\log{\li x \over \varphi(q)}\right). 
\end{equation}
Computations show that formula ({\ref{twologlixphiq}) usually {\it overestimates} $N_{q,r}(x)$. 
At the same time, the right-hand side of ({\ref{twologlixphiq}) is {\it less than} $2\log x$.
Note that a heuristic argument \cite[sect.\,2.3]{kourbatov2017} suggests that 
if the limit $\lim\limits_{x\to\infty}N_{q,r}(x)/\log x$ exists, then the limit is 2.

\subsubsection{Conjecture on the number of records between $x$ and $ex$}\label{conjnum2} 
For any $q>r\ge1$ with $\gcd(q,r)=1$, there exists $\lambda>0$ such that
$$
 N_{q,r}(ex) - N_{q,r}({x}) < (\log\log x)^\lambda  \qquad\mbox{almost always. }
$$
{\em Remark. }
We can also reformulate these conjectures using {\it ``for all $x$ large enough''} 
instead of {\it ``almost always''}.
However, both forms, {\it ``almost always''} and {\it ``for $x$ large enough''}, 
are exceedingly difficult to prove or disprove.

\subsection{Conjecture on $G_{q,r}(x)$ trend. Wolf's conjecture}\label{genwolf}

Our extensive computations suggest the following interesting conjecture
describing the size of maximal gaps between primes in residue classes.

\smallskip\noindent
{\bf Conjecture on the trend of} $G_{q,r}(x)$. \ 
For any $q>r\ge1$ with $\gcd(q,r)=1$, there exist real numbers 
$c_0 = c_0(q)$ and $c_1 = c_1(q)$ such that
\begin{equation}\label{trendconj}
G_{q,r}(x) ~\sim~ T_c(q,x) = 
  {\varphi(q)x\over\li x}\Big(\log{\li x\over\varphi(q)} - c_1\log\log x + c_0\Big) 
  \qquad\mbox{ as }x\to\infty,
\end{equation}
and the difference $G_{q,r}(x)-T_c(q,x)$ changes its sign infinitely often.



\medskip\noindent
{\it Remarks. }

\smallskip\noindent
(i) The trend (\ref{rtrend}) of record {\it random} gaps is a special case of
(\ref{trendconj}), with $c_1=0$.

\smallskip\noindent
(ii) The trend equation (\ref{trend}) can also be regarded as a special case of (\ref{trendconj}),
with $c_1=1$, up to a correction term $O(a)$.

\medskip\noindent
{\bf Wolf's conjecture on maximal gaps between primes }($q=2$). \  
Wolf \cite{wolf1998,wolf2011,wolf2014} gives the following approximation for $G(x)$, 
the {\em maximal prime gap below} $x$, in terms of the prime counting function $\pi(x)$:
%
\begin{equation}\label{wolfconj}
G(x) ~\sim~ g(x) = {x\over\pi(x)}(2\log\pi(x) - \log x + c) \qquad\mbox{ \cite[p.\,11]{wolf2011}}.
\end{equation}
Wolf's reasoning also implies that the difference $G(x)-g(x)$ changes its sign infinitely often.

\medskip 
Conjectures (\ref{trendconj}) and (\ref{wolfconj}) are related. 
Taking $c_1=1$, using the fact that $\varphi(2)=1$ and applying the approximations
$$
{x\over\pi(x)}\approx {x\over\li x}\approx \log x -1, \qquad
\log\pi(x) \approx \log\li x \approx \log x - \log\log x,
$$
for $q=2$ we can rewrite both (\ref{trendconj}) and (\ref{wolfconj})
in this common form:
\begin{equation}
G_{2,1}(x) ~\sim~ \log^2 x - 2\log x\log\log x + O(\log x) \qquad\mbox{ as }x\to\infty.
\end{equation}
So, for $q=2$ and $c_1=1$,
the maximal prime gap predictions (\ref{trendconj}) and (\ref{wolfconj}) agree up to $O(\log x)$. 
It is curious that the {\it logarithm} of twice the twin prime constant $\log(2\Pi_2)$ appears in $c$ of 
(\ref{wolfconj}) \cite[eq.\,31]{wolf2011} --- while it is not immediately clear whether $\Pi_2$ 
might explicitly appear in eqs.\,(\ref{trend}), (\ref{trendconj}). 
Note, however, that there might be multiple meaningful choices of $c_0$; 
e.\,g., three different choices of $c_0$ would map gaps near the corresponding trend $T_c$ to the 
{\it mode}, {\it mean}, or {\it median} of the rescaled gaps $w$. 

\subsection{Consequences of the new conjectures}\label{consequences}

Let us assume for a moment that conjectures of sections \ref{conjnum} and \ref{genwolf} are true.
Assume further that for large $x$ the maximal gaps $G_{q,r}(x)$, after rescaling, 
obey the Gumbel distribution with scale parameter $\alpha\lesssim1$ and mode $\mu\approx0$. 
What consequences do these conjectures entail?

\begin{center}Table~2. \ 
The $c_1$ value in (\ref{trendconj}) allows to predict how often $G_{q,r}(x)>\varphi(q)\log^2x$ \\[0.5em]
\begin{tabular}{ccc}
\hline {\large $\vphantom{1^{1^1}}$}
$c_1$  &  $G_{q,r}(x)<\varphi(q)\log^2x$ &  $G_{q,r}(x)>\varphi(q)\log^2x$ \\
[0.5ex]\hline
     $c_1<-1$ & almost never \vphantom{\fbox{$1^1$}} & almost always                       \\
     $c_1=-1$ & a positive proportion of max.\ gaps  & a positive proportion of max.\ gaps \\
 $-1<c_1\le0$ & almost always                        & almost never but infinitely often   \\
     $c_1>0$  & almost always                        & at most finitely often              \\
\hline
\end{tabular}
\end{center}

Asymptotic analysis (see below) shows that under these assumptions
the $c_1$ value in (\ref{trendconj}) allows us to heuristically predict how often $G_{q,r}(x)>\varphi(q)\log^2x$.
Table 2 lists four particular cases distinguished by the $c_1$ value. 

The observed trend equation (\ref{trend}) suggests that $c_1=1$. 
If our conjectures remain valid for arbitrarily large $x$
{\em and} if there is a limiting Gumbel distribution of record gaps,
this further suggests that there are {\em at most finitely many violations} 
of the ``naive'' Cram\'er conjecture ($G_{2,1}(x)<\log^2x$) and Firoozbakht's conjecture (sect.\,\ref{genfiroozbakht}).
Likewise, for any given $q>r\ge1$, we can expect at most finitely many violations of the generalized
Cram\'er conjecture D in residue classes; see section 
\ref{gencram} for further discussion.

\smallskip\noindent
{\bf Asymptotic analysis.} Let us analyze the most interesting case $c_1>-1$. 
Suppose that $p$ is the endpoint of a record gap $G_{q,r}(p)$,  and $p$ is very large. 
Under the rescaling transformation 
$$
G_{q,r}(x) ~\mapsto~ w = {G_{q,r}(x)-T_c(q,x)\over a(q,x)}
$$
we see from eq.\,(\ref{trendconj}) that 
\begin{itemize}
  \item $G_{q,r}(p)\approx T_c(q,p)$ is mapped to $w\approx0$;
  \item $G_{q,r}(p)\approx \varphi(q)\log^2 p$ is mapped to $w\approx(1+c_1)\log\log p$.
\end{itemize}
The standard Gumbel distribution cdf at $(1+c_1)\log\log p$ is 
$$
\exp(-\exp(-(1+c_1)\log\log p)) = \exp(-1/\log^{1+c_1} p) \approx 1 - 1/\log^{1+c_1} p.
$$
So for record gaps near $p$ we have
\begin{itemize}
  \item $G_{q,r}(p)<\varphi(q)\log^2 p$ with probability $\approx 1-1/\log^{1+c_1} p$;
  \item $G_{q,r}(p)>\varphi(q)\log^2 p$ with probability $\approx 1/\log^{1+c_1} p$.
\end{itemize}
If there was exactly one record gap with endpoint $p\in[e^k,e^{k+1}]$,
then the total number of ``exceptionally large'' record gaps of size $>\varphi(q)\log^2 p$ 
would be estimated as the series
\begin{equation}\label{series1}
\displaystyle\sum\limits_{k=1}^{\infty} {1\over\log^{1+c_1} e^k} = 
\displaystyle\sum\limits_{k=1}^{\infty} {1\over k^{1+c_1}}.
\end{equation}
However, by conjecture \ref{conjnum2} the actual number of record gaps with $p\in[e^k,e^{k+1}]$ 
may be more than one but is below $(\log\log p)^\lambda$ for some $\lambda>0$.
Therefore we can predict that
\begin{equation}\label{series2}
\mbox{ (the total number of records exceeding $\varphi(q)\log^2 p$)} <
\sum\limits_{k=1}^{\infty} {(\log\log e^k)^\lambda\over\log^{1+c_1} e^k} = 
\sum\limits_{k=1}^{\infty} {\log^\lambda k\over k^{1+c_1}}.
\end{equation}
Both series (\ref{series1}), (\ref{series2}) diverge for $c_1\le0$; both series converge for $c_1>0$. 
In particular, if $\lambda$ is a positive integer and $c_1>0$, we have
$$
\sum\limits_{k=1}^{\infty} {\log^\lambda k\over k^{1+c_1}} = |\zeta^{(\lambda)}(1+c_1)|,
$$
where $\zeta^{(\lambda)}$ is the $\lambda$th derivative of the Riemann $\zeta$-function.
In view of the Borel-Cantelli lemma, for $c_1>0$ the convergence of series (\ref{series1}), (\ref{series2}) suggests
that the number of gaps exceeding $\varphi(q)\log^2 p$ is at most finite for any given coprime pair $(q,r)$.

\smallskip\noindent
{\it Remark.} In Cram\'er's probabilistic model of primes \cite[OEIS \seqnum{A235402}]{cramer,kourbatov2014,oeis}
we have a situation corresponding to $c_1=0$; so in Cram\'er's model we should expect that 
$G(p)>\log^2 p$ almost never but still infinitely often.

\section{Generalizations of some familiar conjectures \\ corroborated by experimental data}

In this section we discuss several generalizations of familiar conjectures 
that are related to (and in some cases suggested by) 
the numerical results of Section \ref{numresults}. 
This discussion places our computational experiments in a broader context.
Some of the conjectures (sect.\,\ref{gencram}) have been proposed 
by the author earlier on the {\it PrimePuzzles.net} website \cite{rivera77}.
We always assume that $q>r\ge1$ and $\gcd(q,r)=1$. 

\subsection{Generalized Riemann hypothesis}

The Riemann hypothesis is a statement about non-trivial zeros of the Riemann $\zeta$-function.
An ``elementary'' reformulation of the RH is the {\em prime number theorem} with an
$O(x^{1/2+\varepsilon})$ error term. It states that 
$\pi(x)$, the total number of primes $\le x$, satisfies 
$$
\pi(x) ~=~ {\li x} + O(x^{1/2+\varepsilon}) 
\quad \mbox{ for any } \varepsilon>0, \quad \mbox{ as } x\to\infty.
$$
The generalized Riemann hypothesis is a similar statement concerning non-trivial zeros of 
Dirichlet's $L$-functions. An ``elementary'' reformulation of the GRH says that 
$\pi_{q,r}(x)$, the total number of primes $p=qn+r\le x$, $n\in{\mathbb N}^0$,  satisfies 
$$
\pi_{q,r}(x) ~=~ {\li x\over\varphi(q)} + O_q(x^{1/2+\varepsilon}) 
\quad \mbox{ for any } \varepsilon>0, \quad \mbox{ as } x\to\infty.
$$
Roughly speaking, the GRH means that for large $x$ the numbers $\pi_{q,r}(x)$ and 
$\lfloor\li x/\varphi(q)\rfloor$ almost agree in the left half of their digits.
Thus the GRH justifies our assumption that average gaps between primes $\le x$
in each residue class mod $q$ are nearly the same size, $\varphi(q)x/\li x$. 
Consequently, it is reasonable to expect 
that maximal gaps in each residue class have about the same growth trend and obey 
approximately the same distribution around their trend. 
This is indeed observed in our computational experiments.

(A weaker unconditional result similar to the GRH is the Siegel--Walfisz theorem.
Friedlander and Goldston \cite{friedgold1996} give more results for 
the distribution of primes in residue classes.)

\subsection{Generalizations of Cram\'er's conjecture}\label{gencram}

If $G(x)$ is the maximal gap between {\it primes below} $x$, 
Cram\'er \cite{cramer} conjectured in the 1930s that $G(x)=O(\log^2 x)$.
Clearly $G(x)=G_{2,1}(x)$ for all $x\ge5$.
We give four conjectures (from most to least plausible), each of which can be considered 
a generalization of Cram\'er's conjecture to maximal gaps $G_{q,r}(p)$ between primes 
in a residue class modulo $q$. Everywhere we assume that the maximal gap $G_{q,r}(p)$ ends at the prime $p$.

\subsubsection{``Big-O'' formulation}\label{ocram}

{\bf Generalized Cram\'er conjecture A.} \ 
For any $q>r\ge1$ with $\gcd(q,r)=1$, we have
$$
G_{q,r}(p) ~=~ O(\varphi(q) \log^2 p).
$$
This weakest generalization of Cram\'er's conjecture seems quite likely to be true.

\subsubsection{``Almost-all'' formulation}\label{aacram}
{\bf Generalized Cram\'er conjecture B.} \ 
{\it Almost all} maximal gaps $G_{q,r}(p)$ satisfy
$$
G_{q,r}(p) ~<~ \varphi(q) \log^2 p
\qquad\mbox{ for any } q>r\ge1 \mbox{ with }\gcd(q,r)=1.
$$

For $q=2$ and $r=1$, 
the above two generalizations of Cram\'er's conjecture are compatible with the heuristics of 
Granville \cite{granville}. At the same time, the stronger formulations in 
subsections \ref{limsupcram} and \ref{naivecram} contradict Granville's heuristic reasoning  
which suggests 
$$\limsup_{p\to\infty} {G_{2,1}(p)\over\log^2 p} ~\ge~ 2e^{-\gamma} = 1.12291\ldots
\qquad\mbox{ \cite[p.\,24]{granville}}.
$$

\subsubsection{``Limit superior'' formulation}\label{limsupcram}
{\bf Generalized Cram\'er conjecture C.} \ 
For any $q>r\ge1$ with $\gcd(q,r)=1$, we have
$$
\limsup_{p\to\infty} { G_{q,r}(p) \over \varphi(q) \log^2 p } = 1.
$$

\subsubsection{``Naive'' formulation}\label{naivecram}
{\bf Generalized Cram\'er conjecture D.} \ 
For any $q>r\ge1$ with $\gcd(q,r)=1$, we have
$$
G_{q,r}(p) ~<~ \varphi(q) \log^2 p \qquad\mbox{ if $p$ is large enough.}
$$
Here $p$ is the prime at the end of the maximal gap $G_{q,r}(p)$.

\medskip\noindent
The latter generalization of Cram\'er's conjecture seems a little far-fetched. Nevertheless,
exceptionally large gaps with $G_{q,r}(p)>\varphi(q) \log^2 p$ are extremely rare \cite{rivera77}. 
See {\it Appendix} \ref{appendix-extra-large} for a list of examples of this kind.

\subsection{Generalizations of Shanks conjecture}\label{genshanks}

As above, let $G(p)$ be the maximal prime gap ending at the prime $p$.
Shanks \cite{shanks} conjectured that the infinite sequence of maximal prime gaps  
satisfies the asymptotic equivalence 
$$
G(p)\sim\log^2p \qquad\mbox{ as }p\to\infty \qquad\mbox{ \cite[p.\,648]{shanks}}.
$$
Our numerical experiments suggest the following natural generalizations of the Shanks conjecture
to describe the behavior of $G_{q,r}(p)$.

\subsubsection{``Almost-all'' formulation}
{\bf Generalized Shanks conjecture I.} \ 
For any $q>r\ge1$ with $\gcd(q,r)=1$, there exists 
an infinite sequence $S$ that comprises {\em almost all} maximal gaps $G_{q,r}(p)$
such that every gap in $S$ satisfies the asymptotic equality
$$
G_{q,r}(p) ~\sim~ \varphi(q) \log^2 p \qquad\mbox{ as }p\to\infty.
$$
For $q=2$ and $r=1$, this formulation is compatible with heuristics of Granville \cite{granville}
implying that there should exist an {\it exceptional infinite subsequence} of $G_{2,1}(p)$ satisfying 
$$G_{2,1}(p) ~\sim~ M \log^2p \quad\mbox{ for some } M\ge 2e^{-\gamma}>1.
$$
Indeed, our ``almost-all'' formulation simply means that any exceptional subsequence is very thin;
that is, a zero proportion of maximal gaps $G_{2,1}(p)$ have exceptionally large sizes predicted by Granville.

\subsubsection{Strong formulation}
{\bf Generalized Shanks conjecture II.} \ 
{\em All} maximal gaps $G_{q,r}(p)$ satisfy 
$$
G_{q,r}(p) ~\sim~ \varphi(q) \log^2 p \qquad\mbox{ as }p\to\infty.
$$
This strong formulation contradicts Granville's heuristics cited above.

\subsection{Generalizations of Firoozbakht's conjecture}\label{genfiroozbakht}

Firoozbakht \cite{rivera30} conjectured that $(p_k^{1/k})_{k\in{\mathbb N}}$ 
is a decreasing sequence. Equivalently,
\begin{equation}\label{firoozbakht}
p_{k+1} ~<~ p_{k}^{1+1/k} \qquad\mbox{ for all } k\ge1.
\end{equation}
The conjecture has been verified for all $p_k<10^{19}$ \cite{kourbatov2015v,kourbatov2018w}. 
An independent verification was also performed by Wolf (unpublished); see \cite{sun2013}.
Firoozbakht's conjecture implies
\begin{equation}\label{firoozgap}
p_{k+1}-p_k ~<~ \log^2 p_k - \log p_k - 1  \qquad\mbox{ for all } k>9 
\qquad\mbox{\cite[Th.\,1]{kourbatov2015u}}.
\end{equation}

\subsubsection{The Sun--Firoozbakht conjecture for primes in residue classes}\label{sunfiroozconj}

Z.-W.~Sun \cite[Conjecture 2.3]{sun2013} generalized Firoozbakht's conjecture (\ref{firoozbakht}) as follows:

\medskip\noindent
Let $q > r \ge 1$ be positive integers with $r$ odd, $q$ even and $\gcd(r, q) = 1$. 
Denote by $p_n(r,q)$ the $n$-th prime in the progression $r$, $r+q$, $r+2q\ldots$
Then there exists $n_0(r,q)$ such that the sequence 
$(p_n(r,q)^{1/n})_{n\ge n_0(r,q)}$ is strictly decreasing. 
In particular, one can take $n_0(r,q) = 2$ for $q \le 45$.

\medskip\noindent
{\em Remark.} If the latter conjecture is true, $n_0(r,q)$ may be quite large.
%
%
For example, with $q=486$ and $r=127$ we must have $n_0(r,q)>282866$ because
$$ 
897783067^{1/282866} < 897849649^{1/282867}
$$
(897783067 and 897849649 are the 282866th and 282867th primes $p\equiv127$ mod 486).

Assuming the GRH and reasoning along the lines of \cite{kourbatov2015u}, 
we can prove that the Sun-Firoozbakht conjecture implies
\begin{equation}\label{sunfiroozineq}
{p_{n+1}(r,q)-p_n(r,q)\over \varphi(q)} ~<~ \log^2 p_n(r,q) - \log p_n(r,q) - 1  
\qquad\mbox{ for all $n$ large enough}. 
\end{equation}
This in turn implies the generalized Cram\'er conjectures \ref{ocram}, \ref{aacram}, and \ref{naivecram}, 
as well as a modified form of \ref{limsupcram} with 
$\displaystyle\limsup_{p\to\infty}G_{q,r}(p)/(\varphi(q)\log^2 p)\le1$.

\subsubsection{``Almost all'' form of the Sun--Firoozbakht conjecture}

The Sun--Firoozbakht conjecture of sect.\,\ref{sunfiroozconj} might be excessively strong. 
A weaker, more plausible statement is the ``almost all'' formulation below.

\medskip\noindent
For any coprime $q>r\ge1$, the endpoints $p_{n}(r,q)$, $p_{n+1}(r,q)$
of {\bf\em almost all} record gaps $G_{q,r}(x)$ satisfy inequality (\ref{sunfiroozineq}).

\medskip\noindent
Here, again, $p_n(r,q)$ denotes the $n$-th prime in the arithmetic progression 
$r$, $r+q$, $r+2q\ldots$

\pagebreak
\section{Summary}

The computational experiments described herein took many weeks of computer time.
These experiments allow us to arrive at the following conclusions:
\begin{itemize}
\item The observed growth trend of record (maximal) gaps $G_{q,r}(x)$ is given by eqs.\,(\ref{trend})--(\ref{defb0}). 
      In particular, for record prime gaps ($q=2$) the trend equation reduces to
      $$
      G_{2,1}(x) ~\sim~ \log^2 x - 2\log x\log\log x + O(\log x) \quad\mbox{ as }x\to\infty.
      $$
\item The Gumbel distribution, after proper rescaling, is a possible limit law for $G_{q,r}(p)$.
      The existence of such a limiting distribution is an open question.
\item Almost all maximal gaps $G_{q,r}(p)$ between primes in residue classes modulo $q$ 
      appear to satisfy the asymptotic relation $G_{q,r}(p)\lesssim\varphi(q)\log^2p$.
\item Exceptions with $G_{q,r}(p)>\varphi(q)\log^2p$ are extremely rare (see {\it Appendix} \ref{appendix-extra-large}). 
\item In view of the Borel--Cantelli lemma, our observations suggest that the total number of such exceptions is finite
      for any given coprime pair $(q,r)$ {\em if} there is a limiting Gumbel distribution of rescaled gaps $G_{q,r}(p)$ 
      {\em and if} the trend (\ref{trend}) holds for $x\to\infty$. 
      (In Cram\'er's probabilistic model of primes, on the contrary, we expect infinitely many exceptions 
      $G(p)>\log^2p$.) 
\item We conjecture that the total number $N_{q,r}(x)$ of record gaps up to $x$ is  
      almost always below $\log x (\log\log x)^\lambda$ for some $\lambda>0$.
      It is possible that a stronger conjecture is also true: $N_{q,r}(x)<C\log x$ for some $C>2$.
\end{itemize}
Some of our conjectures may well be undecidable. At any rate, 
such statements appear to be no easier to prove or disprove than the famous Riemann hypothesis.

\pagebreak
\section{Appendix: Details of computational experiments}

Interested readers can reproduce and extend our results using the programs below.

\subsection{PARI/GP program maxgap.gp (ver.\,2.0)}

{\small
\begin{verbatim}
\\ Usage example: q=1000;for(r=1,q-1,if(gcd(q,r)==1,maxgap(q,r,1e12)))

default(realprecision,11) 

\\ li(x) computes the logarithmic integral of x
li(x) = return(real(-eint1(-log(x))))

\\ pmin(q,r) computes the least prime p = qn + r, for n=0,1,2,3,...
pmin(q,r) = forstep(p=r,1e99,q, if(isprime(p), return(p)))

\\ maxgap(q,r,end [,b0,b1,d]) computes maximal gaps g 
\\ between primes p = qn + r, as well as rescaled values (w, u, h):
\\   w and u use Wolf's conjecture (u - without correction term b);
\\   h uses exteme value theory (cf. randomgap.gp).
\\ Results are written on screen and in the c:\wgap folder.
\\ Computation ends when primes exceed the end parameter. 
maxgap(q,r,end,b0=1,b1=4,d=2.7) = {
  re = 0;
  p = pmin(q,r);
  t = eulerphi(q); 
  while(p<end,
    m = p + re;
    p = m + q;
    while(!isprime(p), p+=q);
    while(!isprime(m), m-=q);
    g = p - m;
    if(g>re,
      re=g; Lip=li(p); a=t*p/Lip; Logp=log(p);
      h = g/a-log(Lip/t);
      u = g/a-2*log(Lip/t)+Logp;
      w = g/a-2*log(Lip/t)+Logp-log(t)*(b0+b1/max(2,log(Logp))^d);
      f = ceil(Logp/log(10));
      write("c:\\wgap\\"q"_1e"f".txt", w" "u" "h" "g" "m" "p" q="q" r="r);
                                 print(w" "u" "h" "g" "m" "p" q="q" r="r);
      if(g/t>log(p)^2, write("c:\\wgap\\"q"_1e"f".txt","extra large"));
    )
  )
}
\end{verbatim}
}

\subsection{PARI/GP program randomgap.gp}

{\small
\begin{verbatim}
\\ Usage example: q=1000;for(r=1,q-1,if(gcd(q,r)==1,randomgap(q,r,1e10)))

default(realprecision,11)

\\ exprv(m) returns an exponential random variable with mean m
exprv(m) = return(-m*log(random(1.0)))

\\ randomgap(q,r,end) writes to c:\ygap\ a set of files with record gaps
\\ in a growing sequence of integers p separated by "random" gaps which
\\ are exponentially distributed with mean m = phi(q)*log(p).
\\ The parameter r is included to mimic maxgap(q,r,end).
randomgap(q,r,end) = {
  re=0; p=max(2,r); t=eulerphi(q);
  while(p<end,
    g=ceil(exprv(t*log(p)));
    s=p;
    p+=g;
    if(g>re,
      re=g;
      Lip=real(-eint1(-log(p)));  \\li(p)
      a=t*p/Lip;
      h=g/a-log(Lip/t);
      f=ceil(log(p)/log(10));
      write("c:\\ygap\\rand"q"_1e"f".txt",
            h"  "g" "s" "p" q="q" r="r);
      print(h"  "g" "s" "p" q="q" r="r);
    )
  )
}
\end{verbatim}
}

\subsection{Notes on distribution fitting}

To study the distributions of standardized maximal gaps $G_{q,r}(p)$
we used the distribution fitting software {\tt EasyFit} \cite{easyfit}.
Data files created with {\tt maxgap.gp} or {\tt randomgap.gp}
are easily imported into {\tt EasyFit}: from the {\it File} menu, choose {\it Open},
select the data file, then specify {\it Field Delimiter = space},
click {\it Update}, then {\it OK}.

\medskip\noindent
{\bf Caution:} PARI/GP writes real numbers near zero in a mantissa-exponent format
{\em with a space} preceding the exponent (e.g.\ {\tt1.7874829515 E-5}), 
whereas {\tt EasyFit} expects such numbers {\em without a space} (e.g.\ {\tt1.7874829515E-5}). 
Therefore, before importing into {\tt EasyFit}, search the data files for {\tt" E-"} 
and replace all occurrences with {\tt"E-"}.

\subsection{Exceptionally large gaps: $G_{q,r}(p)>\varphi(q)\log^2p$}\label{appendix-extra-large}

The following table lists exceptionally large maximal gaps $G_{q,r}(p)>\varphi(q)\log^2p$.
No other maximal gaps with this property were found for $p<10^9$, $q<20000$.
Three sections of the table correspond to (i) odd $q,r$; (ii) even $q$; (iii) even $r$.
(Overlap between sections is due to the fact that $\varphi(q)=\varphi(2q)$ for odd $q$.)
No such large gaps exist for $p<10^{10}$, $q<1000$.

\medskip
{\footnotesize
\begin{center}{\small Table~3. \ Exceptionally large maximal gaps: $G_{q,r}(p)>\varphi(q)\log^2p$} \\[0.5em]
\begin{tabular}{rrrrrc}
\hline { \small$\vphantom{1^{1^1}}$}
Gap $G_{q,r}(p)$  &  Start of gap  &  End of gap ($p$) &   $q~~$ & $r~~$ & $G_{q,r}(p)/(\varphi(q)\log^2p$) \\
[0.5ex]\hline
\vphantom{\fbox{$1^1$}} 
(i)~~~~~208650
        &   3415781 & 3624431 & 1605 & 341 & 1.0786589153 \\
316790  &    726611 & 1043401 & 2005 & 801 & 1.0309808771 \\
229350  &   1409633 & 1638983 & 2085 & 173 & 1.0145547849 \\
532602  &    355339 & 887941 & 4227 & 271 & 1.0081862161 \\
984170  &   5357381 & 6341551 & 4279 & 73 & 1.0339720553 \\
1263426 &  10176791 & 11440217 & 4897 & 825 & 1.0056800570 \\
2306938 &  82541821 & 84848759 & 6907 & 3171 & 1.0022590147 \\
3415794 & 376981823 & 380397617 & 8497 & 3921 & 1.0703375544 \\
2266530 & 198565889 & 200832419 & 8785 & 7319 & 1.0335372951 \\ \hline
(ii)~~~~411480
        & 470669167 & 471080647 & 3048 & 55 & 1.0235488825 \\
208650  &   3415781 & 3624431 & 3210 & 341 & 1.0786589153 \\
316790  &    726611 & 1043401 & 4010 & 801 & 1.0309808771 \\
229350  &   1409633 & 1638983 & 4170 & 173 & 1.0145547849 \\
657504  & 896016139 & 896673643 & 4566 & 2563 & 1.0179389550 \\
1530912 & 728869417 & 730400329 & 6896 & 3593 & 1.0684247390 \\
532602  &    355339 & 887941 & 8454 & 271 & 1.0081862161 \\
984170  &   5357381 & 6341551 & 8558 & 73 & 1.0339720553 \\
1263426 &  10176791 & 11440217 & 9794 & 825 & 1.0056800570 \\
2119706 & 665152001 & 667271707 & 10046 & 6341 & 1.0223668231 \\
1885228 & 163504573 & 165389801 & 10532 & 5805 & 1.0000704209 \\
1594416 & 145465687 & 147060103 & 13512 & 9007 & 1.0026889378 \\
2306938 &  82541821 & 84848759 & 13814 & 3171 & 1.0022590147 \\
3108778 & 524646211 & 527754989 & 15622 & 12585 & 1.0098218219 \\
1896608 &    164663 & 2061271 & 16934 & 12257 & 1.0598397341 \\
3415794 & 376981823 & 380397617 & 16994 & 3921 & 1.0703375544 \\
2266530 & 198565889 & 200832419 & 17570 & 7319 & 1.0335372951 \\
2937868 &  71725099 & 74662967 & 17698 & 12803 & 1.0103309882 \\
2823288 &  37906669 & 40729957 & 18098 & 9457 & 1.0162761199 \\
2453760 &  11626561 & 14080321 & 18176 & 12097 & 1.0107626289 \\
3906628 & 190071823 & 193978451 & 18692 & 11567 & 1.1480589845 \\ \hline
(iii)~~~657504
        & 896016139 & 896673643 & 2283 & 280 & 1.0179389550 \\
2119706 & 665152001 & 667271707 & 5023 & 1318 & 1.0223668231 \\
3108778 & 524646211 & 527754989 & 7811 & 4774 & 1.0098218219 \\
1896608 &    164663 & 2061271 & 8467 & 3790 & 1.0598397341 \\
2937868 &  71725099 & 74662967 & 8849 & 3954 & 1.0103309882 \\
2823288 &  37906669 & 40729957 & 9049 & 408 & 1.0162761199 \\
3422630 &    735473 & 4158103 & 14881 & 6304 & 1.0368176014 \\
3758772 & 144803717 & 148562489 & 15927 & 11360 & 1.0000152764 \\
3002682 &   8462609 & 11465291 & 16869 & 11240 & 1.0107025944 \\
8083028 & 344107541 & 352190569 & 19619 & 9900 & 1.1134625422 \\
\hline
\end{tabular}
\end{center}
}

\section{Acknowledgments}
I am grateful to all contributors and editors
of the websites {\it OEIS.org} and {\it PrimePuzzles.net}. Thanks also to
Marek Wolf and Zhi-Wei Sun for putting forward many interesting conjectures 
\cite{sun2013}--\cite{wolf2016}.
I am especially grateful to Prof.\,Wolf for his very useful comments via email.


{\small

\bigskip
\hrule
\bigskip

\noindent 2010 {\it Mathematics Subject Classification}: 11N05.

\noindent \emph{Keywords: }
Cram\'er conjecture, Firoozbakht conjecture, Gumbel distribution, 
prime gap, residue class, Shanks conjecture, Wolf conjecture.

\bigskip
\hrule
\bigskip

\noindent (Concerned with sequences
 \seqnum{A005250},
 \seqnum{A084162},
 \seqnum{A235402},
 \seqnum{A235492},
 \seqnum{A268799},
 \seqnum{A268925},
 \seqnum{A268928},
 \seqnum{A268984},
 \seqnum{A269234},
 \seqnum{A269238},
 \seqnum{A269261},
 \seqnum{A269420},
 \seqnum{A269424},
 \seqnum{A269513},
 \seqnum{A269519}.)

\bigskip
\hrule
\bigskip
\end{document}